\documentclass[graybox]{svmult}

\usepackage{mathptmx}       
\usepackage{helvet}         
\usepackage{courier}        
\usepackage{type1cm}        
\usepackage{makeidx}         
\usepackage{graphicx}        
\usepackage{multicol}        
\usepackage[bottom]{footmisc}
\usepackage{rotating}

\newtheorem{thm}{Theorem}

 \newtheorem{prop}{Proposition}
 \newtheorem{defn}{Definition}

\makeindex             

\begin{document}

\title*{A Copula Approach to Inventory Pooling Problems with Newsvendor Products}
\author{Burcu Ayd{\i}n and Kemal Guler and Enis Kay{\i}\c{s}}
\institute{Burcu Ayd{\i}n \at HP Labs, 1501 Page Mill Rd Palo Alto CA 94304, \email{aydin@hp.com}
\and Kemal Guler \at HP Labs, 1501 Page Mill Rd Palo Alto CA 94304, \email{kemal.guler@hp.com}
\and Enis Kay{\i}\c{s} \at HP Labs, 1501 Page Mill Rd Palo Alto CA 94304, \email{enis.kayis@hp.com}}
%
%
\maketitle

\abstract*{This study focuses on the inventory pooling problem under the newsvendor framework. The specific focus is the change in inventory levels when product inventories are pooled. We provide analytical conditions under which an increase (or decrease) in the total inventory levels should be expected. We introduce the copula framework to model a wide range of dependence structures between pooled demands, and provide a numerical study that gives valuable insights into the effect of marginal demand distributions and dependence structure on the effect of pooling to inventory levels.}

\abstract{This study focuses on the inventory pooling problem under the newsvendor framework. The specific focus is the change in inventory levels when product inventories are pooled. We provide analytical conditions under which an increase (or decrease) in the total inventory levels should be expected. We introduce the copula framework to model a wide range of dependence structures between pooled demands, and provide a numerical study that gives valuable insights into the effect of marginal demand distributions and dependence structure on the effect of pooling to inventory levels.}

\section{Introduction: Inventory Pooling Problem}

In this paper, we study the inventory pooling problem using the classic newsvendor framework. The newsvendor problem occurs when for a given item,
the inventory level is decided before the realization of its demand. Therefore the optimal inventory level needs to be decided based on the distribution of the stochastic demand $D$. Typically, unsold items at the end of the period are assumed to be either discarded or salvaged. The solution of the newsvendor problem is well-known: a quantile of the demand distribution depending on price and cost of the item is the stock level that is optimal in terms of profit.

The pooling problem occurs when the decision makers have the option to combine inventories for an item that serves multiple demand sources. The pooling could be in the form of determining one physical inventory holding location that will serve multiple locations, setting up a quick transshipment modes between different inventory locations (therefore allowing to plan inventories together), or even designing two products so that they are substitutable for each other if need arises. This effort has a clear reward: It is a well-known fact in the literature that pooling always leads to higher profits. (see, for example, Corbett and Rajaram (2006)). However, the optimal inventory levels in the system may increase or decrease after pooling.

Our paper primarily focuses on the change in optimal inventory levels when demands from multiple sources are pooled. The change in inventory levels is an important decision factor. Contrary to common intuition, pooling may result in a decrease or increase in total inventory levels. The pooling decision may bring additional costs that depend on the targeted inventory level. The costs could be due to adjusting warehouse capacities, re-design costs, etc., and profits may include reduced stock-out rates therefore higher customer satisfaction. These should be carefully weighed together with the profit increase due to pooling. Furthermore, after the pooling decision is made, adjusting the inventory levels to the new optimal levels is important in achieving the higher profits. The new optimal levels depend on the demand for the product in each channel, and how these channels affect each other. We investigate how pooled inventory levels are affected by marginal distributions of product demands and the dependence structure between them.

The question we tackle in this paper is mentioned, though not solved, by Corbett and Rajaram (2006):

\begin{quotation}
Most of this literature in inventory pooling, ... , focuses on the impact of pooling on expected profits. A related, but usually more intractable problem, concerns the effect of pooling on optimal inventory levels. We do not consider that question here, though some work, including Eppen (1979), Erkip et al. (1990), and Van Mieghem and Rudi (2002) do address that issue under more restrictive distributional assumptions than ours. So far, the work related to pooling of inventories has generally lacked a formal mechanism for assessing the impact of dependence on the value of pooling when demands are nonnormal. Whenever dependence has been explicitly included, it has generally been in the context of bivariate or multivariate normal demands.
\end{quotation}

Other previous studies on the subject handled certain cases where the demand distributions of the channels and their relationship can be explained by well known multi-variate distributions. Often independent and identically distributed (IID) demand is assumed. For example, Gerchak and Mossman (1992) assume IID exponentially distributed marginals, and Yang and Shrage (2009) studies IID right-skewed marginals. This approach provides mathematical tractability. In real life applications, however, product demands are neither identical nor independent. In this paper we take up this problem and show that the theory of copulas provides a powerful, tractable yet rigorous, framework to address the effect of relaxing both independent and identical demand assumptions on the optimal pooled inventory levels. They also allow us to analyze a very wide range of different dependence structures that may not fit into any of the well-known multi-variate distributions.

The details of the model we consider are as follows. The cost of stocking each unit is $c$. For each demand unit that can be satisfied from inventory, a revenue of $p$ is made. Unsatisfied demand is lost as well as the overstocked items. The objective is to decide the inventory level $Q$ that will maximize the expected total profit. It is well known that the optimal inventory level is a quantile of the demand distribution, i.e.:
\begin{equation}
    F^{-1}\left(\frac{p-c}{p}\right)=arg \max_{Q} \left\{ p E_{D} \left[min(D,Q)\right]-c Q \right\},
\end{equation}
where $F(.)$ is the distribution function of demand.

In the stylized inventory pooling problem, two identical items with uncertain demands, $D_1$ and $D_2$ are considered. These items have the same unit profit and unit stocking cost. The decision maker has two options: Keeping dedicated inventory to satisfy the demand of each item, or holding a single inventory for the aggregate demand, $D_1+D_2$. It has been shown that pooling is a better option, however one still needs to decide on the optimal inventory levels. In the first option, the optimal inventory in the system can be shown to be $F_1^{-1}(t)+F_2^{-1}(t)$, where $F_i(.)$ is the marginal distribution function of $D_i$ and $t:=\frac{p-c}{p}$ is defined as the margin ratio. It is easy to see that this quantity is independent of the dependence structure between the demands. On the other hand, the optimal pooled inventory level, $F_{1+2}^{-1}(t)$, depends not only on the marginal demand distributions, but also on the dependence structure between $D_1$ and $D_2$ (where $F_{1+2}(x):=Pr(D_1+D_2\leq x)$).

From a practical point of view, the manager knows pooling is a better option, but he needs to decide whether to keep more or less total inventory as a result of that decision. If pooling requires higher levels of total inventory, we say that \emph{pooling effect is positive}. Similarly, \emph{pooling effect is negative} when pooled inventory level is lower than the dedicated inventory. In other words, we define the pooling effect as $F_{1+2}^{-1}(t)-F_1^{-1}(t)-F_2^{-1}(t)$.

\section{Literature Review}

The inventory pooling problem has been studied extensively in the operations management literature. For many of these studies, the main focus has been the profit comparison under various settings. A smaller number of studies take up the problem of determining the pooled inventory levels.

The earliest and most well-known reference on the pooling problem is Eppen (1979). This study considers the pooling problem when product demands are jointly distributed with multi-variate normal distribution with a known covariance matrix. He shows that the centralized system brings cost savings, and the magnitude of these savings depend on the correlation: the lower the dependence, the higher the savings.

Since Eppen (1979), costs and benefits of inventory pooling are investigated under various other settings. See Gerchak and He (2003) and Alfaro and Corbett (2003) for recent reviews.

Netessine and Rudi (2003) focus on the inventory centralization problem for substitutable products. Substitution is the technical equivalent of pooling when full substitution without stock out penalties are allowed. They show that, when substitution is allowed, it is possible that the optimal inventory levels may increase for some items that are being pooled. However, they give results on the levels of individual items, and they do not provide any result on the total inventory level of the items being pooled under centralization.

Erkip et al. (1990) takes up a similar question: the centralization of inventory ordering policies under the newsvendor framework. They investigate the effect of correlation between normally distributed demands of items that can be centralized. They conclude that \textquotedblleft the effect of correlation can be highly significant, resulting in significantly larger amounts of safety stock for optimal control compared to the no-correlation case."

In their 2003 paper, Gerchak and He investigate the effect of demand variances on the pooling savings. They provide a framework in which an increase in demand variability always increases the savings achieved by combining these demands. They do not require the demand distributions to be independent for their result to hold. However, they do not study how the combined inventory levels are affected by variability.

Alfaro and Corbett (2003) ask an interesting question: if the inventory levels are not optimal in a current setting, would pooling still bring savings? They investigate the profits coming from pooling under non-optimal inventory levels, and compare this to the benefit of optimizing the separate inventory levels rather than pooling them. They find conditions under which it is better to optimize the inventory levels of dedicated setting, and conditions where pooling only will be more profitable.

In inventory pooling literature the effect of dependence on the optimal inventory levels has been studied assuming multi-variate normal demands. The only exception is Corbett and Rajaram (2006). Corbett and Rajaram (2006) use copulas to model the dependence structure between demands. As noted by the authors, they focus on the impact of pooling on expected profits. This focus is particularly essential as the results of superiority of pooling rely critically on the ability to find the optimal inventory levels.

The small number of studies that focus on the pooled inventory levels provide examples in which pooling leads to higher inventory levels contrary to the earlier intuitions. For example, Pasternack and Drezner (1991) show that this comparison depends on the transfer revenue.  Transfer revenue is the profit that comes from the substitution of one product when the other's inventory is depleted. Their cost structure for the substituted amount is different than the original costs, therefore their results are not directly comparable to the studies where pooling is understood as full-substitution, where costs do not change if parts are substituted.

Gerchak and Mossman (1992) conclude that, contrary to the prevalent intuition, pooling may lead to higher inventory levels when demands are independent and identically distributed with an exponential distribution and price per unit and the ratio of cost of underage to cost of overage is sufficiently low. They show this using a numerical counter-example where the demands have exponential distribution. However, they do not provide any generalized findings in terms of providing conditions or distributions which would imply higher pooled inventory levels.

In a recent paper, Yang and Shrage (2009) define the case in which pooling increases inventory levels as \textquotedblleft inventory anomaly". They focus on IID right skewed demand distributions as marginals. They claim that for any two IID right-skewed demand distributions, there exists a range of the margin ratio $\frac{p-c}{p}$ where pooling leads to higher inventory levels. Moreover, for any newsvendor ratio $\frac{p-c}{p}\geq 0.5$, a right-skewed distribution of IID marginals that leads to higher pooled inventory levels exist. Their result is important in the sense that they describe certain conditions where the \textquotedblleft inventory anomaly" can be expected.

Our paper attempts to provide a much more general framework, where the level of pooled inventory can be found under any demand distribution and dependence structure through the use of copulas. Our numerical analysis provides examples of some well-known copulas and marginal distributions that can be used.

\section{Comparison of Inventory Levels}

In this section, we explore a qualitative question: How does the sign of pooling effect change as the margin ratio $t$ varies? The following Proposition sheds light into this question:

\begin{prop}
\label{PropIFONLYIF}
Let $P(t)=F_{1+2}^{-1}(t)-F_1^{-1}(t)-F_2^{-1}(t)$. Assume that $P(t)$ has a unique root $t_0$ in $(0,1)$. Then, $t_0$ is a threshold such that the pooling effect is negative at $t$ if and only if $t>t_0$.
\end{prop}

Proof of Proposition \ref{PropIFONLYIF} follows from Theorem $1$ of Liu and David (1989). The proposition stipulates that if $P(t)$ has a unique root in $(0,1)$, then pooling effect can only change sign from negative to positive as $t$ increases, the critical threshold being $t_0$.
This threshold depends on both the marginals and the dependence structure of joint demand distribution. However, we can characterize this value under certain settings. First, it is easy to verify that for multi-normal family of demand distributions, this critical threshold is always 0.5. That is, regardless of specific parameters that describe a multi-normal demand, pooling leads to a lower inventory level if and only if margin ratio is higher than 0.5. One can extend this \textquotedblleft detail-free" threshold result to other distributions from the same family.

\begin{prop}
\label{CorElliptic}
If $(D_1,D_2)$ follows a distribution in the elliptical family\footnote{The elliptical family includes well-known distributions such as Normal, Laplace, Student-t, Cauchy, and Logistic among others. }, then pooling leads to lower inventory if and only if the margin ratio is higher than $0.5$.
\end{prop}

The proof of this proposition follows trivially from Theorem $6.8$ of McNeil et al. (2005).


For the next result, we use the following definition:
\begin{definition}
A distribution function $F$ is regularly varying at minus infinity with tail index $\alpha>0$ if,
\[\lim_{t\rightarrow \infty} \frac{F(-tx)}{F(-t)}=x^{-\alpha}\ \ \ \  \forall x>0\]
\end{definition}

The next proposition shows that if the distribution of demand is regularly varying, then the sign of pooling effect depends on the tail properties of the demand distribution.
\begin{prop}
\label{deneme2}
Assume that the tail probability of the joint demand distribution to be negligible compared to those of marginal demand distributions. Moreover, let $D_1$ and $D_2$ are identically distributed with regularly varying distribution functions with the same tail index $\alpha$. There exists a threshold $0<t_0<1$ such that if the margin ratio $t$ is greater than or equal to $t_0$ then:
\begin{itemize}
  \item the pooling effect is negative if $\alpha >1$.
  \item the pooling effect is positive if $0<\alpha <1$.
\end{itemize}
\end{prop}
The proof is from Theorem $10$ of Jang and Jho (2007).

It is possible to have $t_0=0$ which implies that pooling leads to higher inventory for all margin ratios. For example, when demands are independent and identically distributed with Pareto distributions, which has infinite mean, then the threshold becomes $0$.

Having established that pooling may lead to either higher or lower inventory levels and some conditions for positive an negative pooling effect, we next investigate the effect of characteristics of the demand uncertainty on inventory levels. With multiple products, we need to study the effect of the marginal demand distributions, as well as the dependence structure between these demands. Towards this end, copula representation provides a unified and rigorous approach which we provide a short overview next.

\section{Brief Overview of Copula Theory}

The joint distribution of demand is critical in understanding the behavior of optimal pooled inventory levels. There are two components of a joint distribution: the marginal distributions of each demand source, and the dependence between these demand sources. In order to study the effect of these components independently, we introduce the copula theory. \footnote{Nelsen (1999) provides an excellent general introduction to the theory of copulas.} Copulas join the univariate marginal distributions of individual random variables to arrive at the joint distribution function for these variables.

\begin{defn}
A d-dimensional copula $C(u_1,\ldots,u_d)$ is a distribution function on $[0,1]^d$ with standard uniform marginal distributions.
\end{defn}

McNeil et al. (2005) shows that a function $C$ with the following properties is a copula:
\begin{enumerate}
  \item $C(u_1,\ldots,u_d)$ is increasing in each component $u_i$.
  \item $C(1,\ldots,1,u_i,1,\ldots,1)=u_i$ for all $i \in \{1,\ldots,d\}$, $u_i \in [0,1]$.
  \item For all $(a_1,\ldots,a_d)$, $(b_1,\ldots,b_d)$ in $[0,1]^d$ with $a_i \leq b_i$, we have:
  \[ \sum_{i_1=1}^{2} \cdots \sum_{i_d=1}^{2} (-1)^{i_1+\cdots+i_d}C(u_{1i_1},\ldots,u_{di_d})\geq 0\]
  where $u_{j1}=a_j$ and $u_{j2}=b_j$ for all $j \in \{1,\ldots,d\}$.
\end{enumerate}

Sklar's Theorem shows that when the marginal distributions are continuous, then the copula is unique.

\begin{thm}[Sklar's Theorem]
\label{TheoremSklar}
Let $F(x_1,x_2,...,x_n)$ be an n-dimensional joint distribution with continuous marginals $F_1(x_1),...,F_n(x_n)$. Then the joint distribution has a unique copula representation given by
\begin{equation}
F(x_1,x_2,...,x_n)=C(F_1(x_1),...,F_n(x_n))
\end{equation}
\end{thm}

Sklar's Theorem provides a powerful technique that enables the separation of marginal distributions from the dependence structure. Since one can fix or
vary the marginals and the copula separately, a rich class of stochastic models can be constructed. Copula-marginal representation of the joint distribution of a set of random variables has been used in a variety of application areas from decision and risk analysis to finance.

Some of the most commonly used examples of the copula include the product, Gaussian, and Archimedean family copulas. The product copula models the independent marginals case. The most commonly used multi-variate distribution, Normal, can be uniquely represented by normal marginals and a Gaussian copula.

Three important copulas within the Archimedean family are Gumbel, Clayton and Frank. An important property that can be modeled using some Archimedean family copulas is the asymmetry around the mode. With these copulas, the dependence structure varies along the different section of the distribution tails. Two copulas that have this property are Gumbel and Clayton. Gumbel distribution is given by:
\[
C_{Gu}(u_1,...,u_n)=e^{-((-ln(u_1))^\theta+...+(-ln(u_n))^\theta)^{1/\theta}},
\]
and Clayton copula is given by:
\[
C_C(u_1,...,u_n)=\left(1-n+\sum_{i=1}^{n}u_i^{-1/\theta}\right)^{-\theta}.
\]
Gumbel copula could be used when the dependence is higher in the right tail, and Clayton could be used when the dependence is higher in the left tail. Clayton copula exhibits comparatively greater dependence in the left tail. Finally, the Frank copula is given as:
\[
C_F(u_1,...,u_n)=\log_{\alpha}\left(\frac{\prod_{i=1}^{n}\alpha^{u_i}-1}{(\alpha-1)^{(n-1)}}+1\right)
\]
Frank copula is symmetric; it exhibits dependence on both tails. In our numerical analysis, we will focus on these three Archimedean copulas. The most important shape qualities of the popular Gaussian copula is already carried by the Frank copula, moreover, Frank copula can represent negative dependence structures as well.

While all of these copula functions represent different two-dimensional structures, they can be compared through a summary scalar of dependence. One such scalar used commonly is called Kendall's $\tau$. For a joint distribution, Kendall's $\tau$ is independent of the marginals and only depends on the copula. It varies between $[-1,1]$, while $-1$ represents perfect negative correlation, $0$ represents lack of correlation and $1$ represents perfect positive correlation. Given a two-dimensional copula function, the associated Kendall's $\tau$ can be found through the following formula:
\begin{equation}\label{kendall}
\tau=4\int_0^1\int_0^1 C(u,v)dC(u,v) -1
\end{equation}
This formula can be found in Kaas et al. (2009).

Other scalar measures of dependence also exist. Out of those, we do not use Pearson's $r$, since it only measures linear dependence and is not a robust measure of nonlinear dependence cases. Spearman's $\rho$ and Blomqvist's $\beta$ are two others that are also common and can measure nonlinear dependence. These could have been used instead of Kendall's $\tau$. However our numerical results would have come out very similar, therefore we limited our analysis to Kendall's $\tau$ only. A wide discussion of these measures can be found in McNeil et al. (2005) and Nelsen (1999).

\section{Numerical Analysis}

In this section, we independently study the effect of identical versus nonidentical and symmetric versus asymmetric marginal demand uncertainties, as well as different types of copulas with varying forms and levels of tail dependence. Our focus is to investigate whether a unique threshold exist beyond which pooling leads to higher inventory levels, the value of this threshold, and the magnitude of the pooling effect. We connect our observations to managerial insights and complement the existing work on optimal pooled inventory levels.

The first factor in determining the inventory levels is the marginal distribution of each demand source. To understand the effect of skewness in marginal demand sources, we use the beta family. The support for the standard beta family is $[0,1]$, hence the optimal total inventory level is always between $0$ and $2$. To study the effect of left or right skewness in the marginal demand, we use $\beta(2,8)$ and $\beta(8,2)$ respectively; we keep the variance fixed by only exchanging the two parameters of the distribution. The case with equal parameters ($\beta(5,5)$) represents the distribution example where the density is symmetric at both tails. Many sources use Normal distribution for this purpose which cannot model skewed cases. Another setting we will investigate is
when the marginals are not identical, by keeping one marginal fixed while varying the other.

The knowledge of marginal distributions is sufficient to determine the optimal dedicated inventory level. However the optimal pooled inventory level also depends on the dependence structure. To understand the effect of dependence structure independent of the level of dependence, we fixed Kendall's $\tau$ to $4$ different values ($0$, $0.2$, $0.5$, and $0.8$) and computed the corresponding copula parameters using Equation \ref{kendall} for the copula families used. For Frank copula that allows negative correlation, we followed the same procedure on the negative side.

Copulas are very powerful tools for representing a wide range of dependence structures. Figure 1 depicts some main copula functions that we will be using in our analysis: Gumbel, Clayton and Frank. All of the these can model strong dependence as well as weak dependence.

\begin{figure}[ptb]\label{copulaspic}
\begin{center}
\includegraphics[scale=0.45]{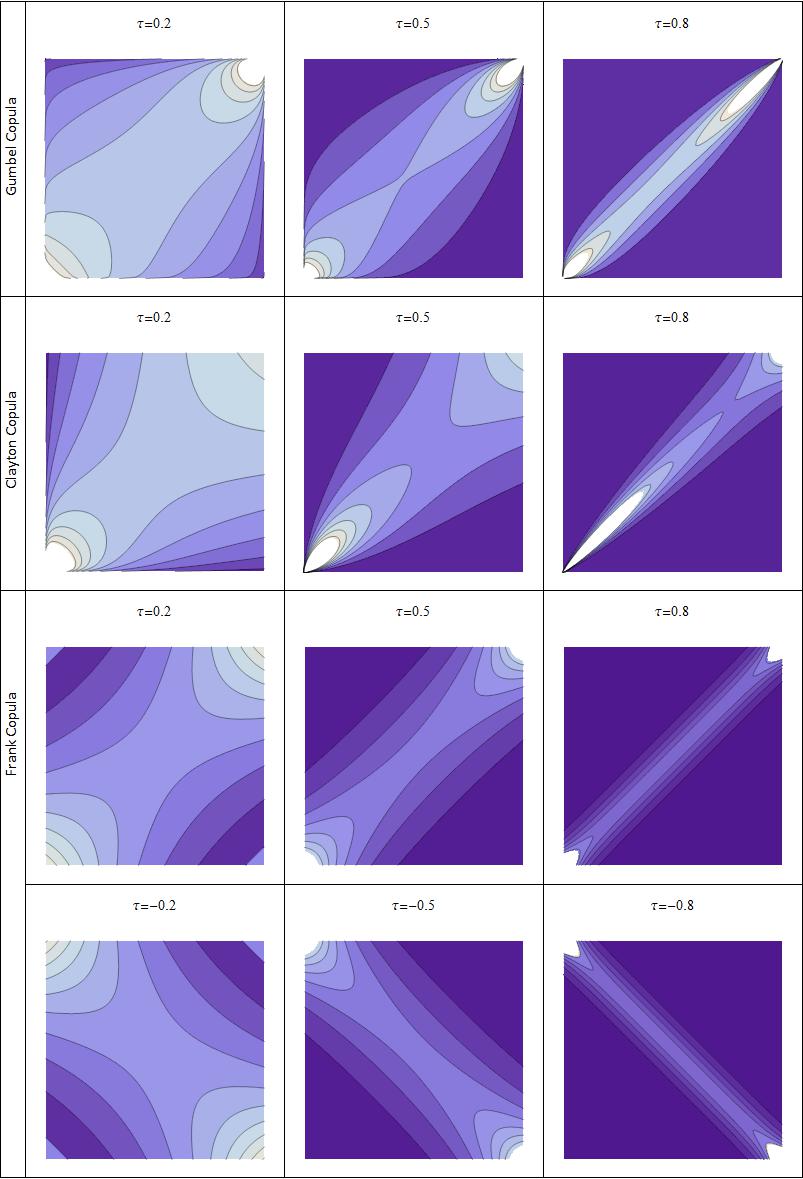}
\end{center}
\caption{Example dependence structures between two dependent random variables represented as contour plots of their copulas. Three Archimedean copula functions, Gumbel, Clayton and Frank are shown. Lighter shades represent the areas with higher density, and darker shades represent areas with lower density. The function parameters are selected such that the copulas depict the dependence structures under Kendall's $\tau=0.2$, $0.5$ and $0.8$ for positive dependence cases (first three rows), and Kendall's $\tau=-0.2$, $-0.5$ and $-0.8$ for negative dependence with Frank copula (fourth row).}
\end{figure}

Comparing the graphs of different copulas under the same Kendall's $\tau$ in Figure 1, we can see that similar levels of \textquotedblleft correlatedness" can exist in very different dependence structures. As Kendall's $\tau$ increases, the densities tend to concentrate around the $45$ degree line. Gumbel copula is appropriate to model cases in which it is slightly more likely that high-level demands are correlated (i.e, higher the dependence on the right top quadrant). Clayton copula models cases where low-level demands are more correlated, perhaps due to unfavorable market conditions that affect all demand sources (i.e, higher the dependence on the left bottom quadrant). Frank copula, on the other hand, shows a more dispersed structure and model cases where dependence is similar in high and low level demands (i.e., it is symmetrical at both tails).

The combined affect of marginals and copulas is what drives the magnitude and sign of the pooling effect at any margin ratio. We will not give the joint density plots of all the marginal-copula pairs that will be used in the
numerical analysis, but for illustrative purposes, the density plots belonging to Gumbel copula are given in Figure 2. Other density plots reveal similar observations, so they will be omitted to save space.
\begin{figure}[ptb]\label{joint}
\begin{center}
\includegraphics[scale=0.7]{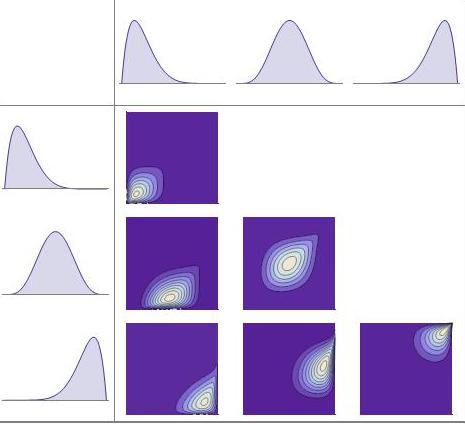}
\includegraphics[scale=0.7]{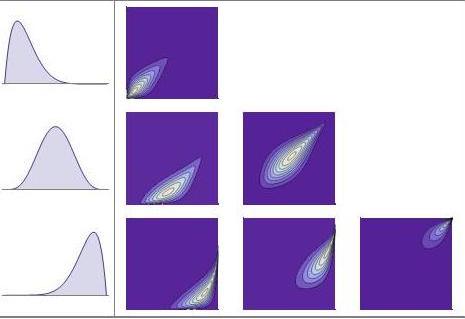}
\includegraphics[scale=0.7]{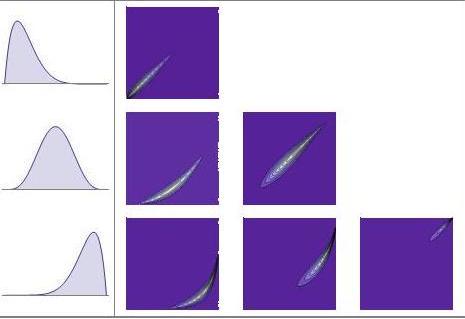}
\end{center}
\caption{The density functions of joint distributions obtained by combining Gumbel copula at three different Kendall's $\tau$ levels, and three different marginals ($\beta(2,8)$,$\beta(5,5)$ and $\beta(8,2)$. On top row and the leftmost column, the plots of marginal densities used are given. The top panel shows six different combinations of these marginals under the Gumbel copula with Kendall's $\tau=0.2$. Middle panel shows the same combinations under Gumbel copula with $\tau=0.5$, and bottom panel is $\tau=0.8$. The effect of higher dependence can be seen from top to bottom: higher dependence concentrates the mass around its center, and reduces dispersion. The effect of non-symmetry of Gumbel copula also becomes visible.}
\end{figure}

To illustrate how both the pooled inventory level and the sum of dedicated inventory levels change under margin ratio, we present Figure 3. This graph gives the intuition on how the dedicated and also pooled inventory levels change when margin ratio changes; i.e., when risk taking is more or less costly.
We see that while the total dedicated inventory level steadily rises with respect to the margin ratio, pooled inventory level is more robust when margin ratio is medium and more sensitive when margin ratio is either too small or too high.
The effects seen on this graph are consistent with the behaviors we observe in Figures 4, 5, 6 and 8.
Similar graphs under different copulas and marginals contain similar trends with respect to the shape of inventory level curves, so they will not be presented here.
\begin{figure}
\begin{center}
\includegraphics[scale=0.7]{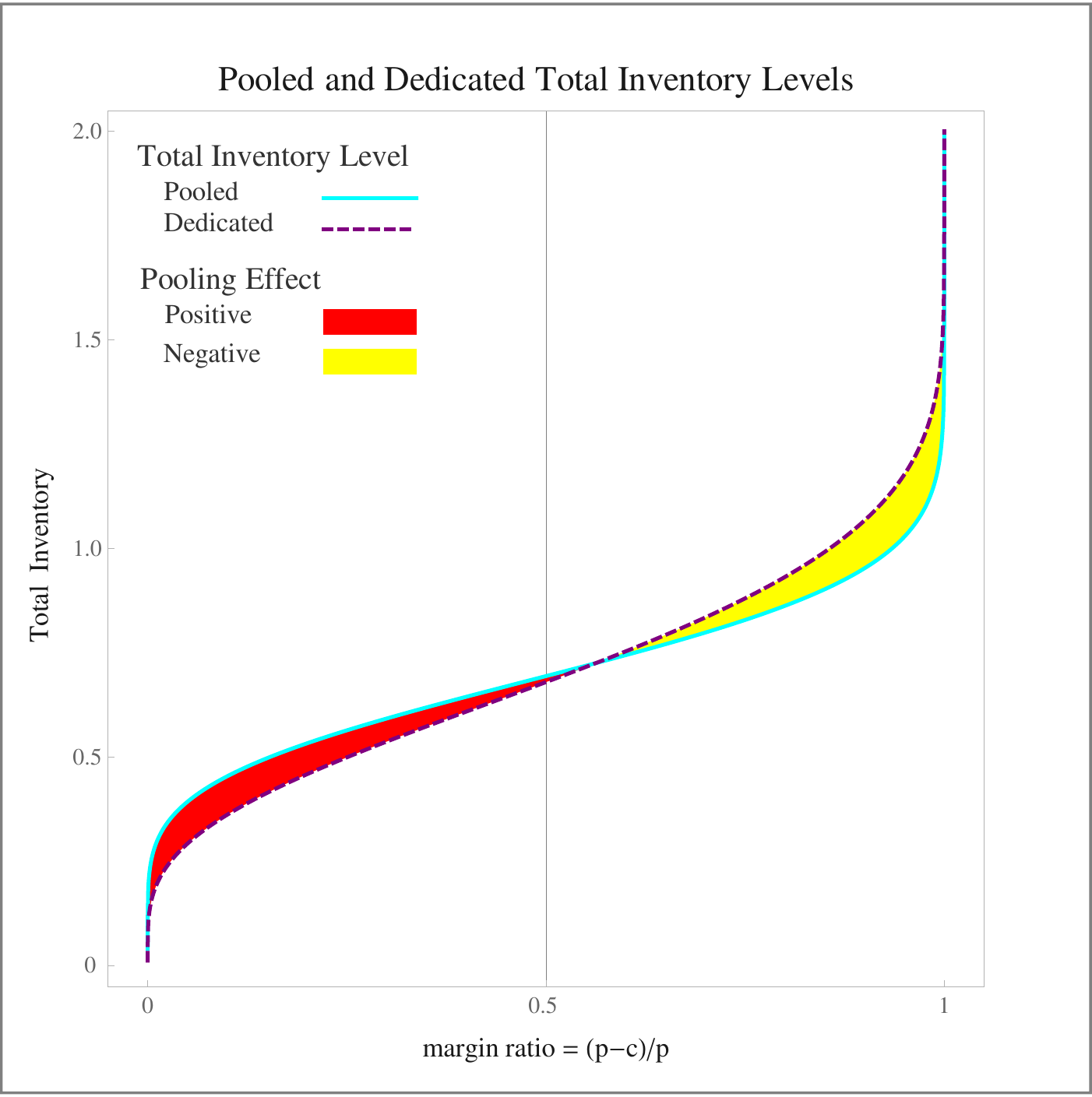}
\end{center}
\caption{The graph of dedicated (dashed purple line) and pooled inventory levels (solid blue line) where the X axis represents the margin ratio. Product demands are identically distributed with Beta where $\alpha=2$ and $\beta=4$. Their dependence structure is given by Frank copula with parameter $\alpha=100$ (Kendall's $\tau=-0.43$ for this parameter). The X axis changes from $0$ to $1$ where $0.5$ threshold is marked with a vertical gray line.}
\end{figure}

\subsection{Effect of Marginal Demand Distribution}
The skewness of the marginal demand distribution affect both the threshold and the magnitude of pooling effect. First, the threshold $t_0$ may be either less than or greater than the critical value $0.5$ as the skewness change. In general, left skewness in any marginal tend to increase the threshold, while right skewness have the opposite effect. This means that when a demand source is more likely to be low than high, then pooling leads to higher inventory for even smaller levels of the margin ratio. This follows from the fact that left skewed distributions concentrate the mass to lower left area of the joint distribution density, while right skewed ones concentrate on the upper right area.

For two identical right-skewed marginals, it is clear from the plots that pooling effect is positive for small values of $t$. This is consistent with the result of Yang and Schrage (2009) on the existence of positive pooling effect for two IID right skewed distributions with small $t$. Our numerical results extend this finding to non-identical marginal demands. We find that as one of the marginals change from left skewed to right skewed while keeping the other one fixed, the threshold decreases, implying a positive pooling effect over a larger set of margin ratio. We also observe that the left skewness of marginals increases the threshold. Finally, we find that the case with one demand marginal being left-skewed and the other being right-skewed leads to similar results to the case wherein both marginals are symmetric around the mean.

The magnitude of the pooling effect decreases as marginals change from being left skewed to right skewed, for any or both of the marginals. When both of the marginals are left skewed we observe that pooling changes the inventory levels most, especially when margin ratio is small.

\subsection{Effect of Dependence Structure}
A comparison across the different copulas gives us interesting insights into seeing the effect of dependence structure, apart from the level of dependence itself, on the inventory levels.

We start with some general conclusions that can be drawn from Figures 4, 5, 6 and 8. First, a stronger dependence measured by a high Kendall's $\tau$ leads to a pooling effect that is smaller in absolute value. This is expected, as we know that for co-monotone demand distributions, the sum of quantiles of individual demand distributions is equal to the quantile of the sum of the two demand distributions, implying no pooling effect.

\subparagraph{Gumbel copula} For the high-dependence case ($\tau=0.8$), we see that pooling does not have a strong effect on inventory levels for most margin ratios. This is consistent with previous knowledge. One exception is when the margin ratio is very low: for any marginal density combination, low margin ratios result in positive pooling effect. Another observation is that as the dependence decreases, the threshold value increases, implying a positive pooling effect for even smaller values of margin ratio when dependence is high.
\begin{figure}
\begin{center}
\includegraphics[scale=0.4]{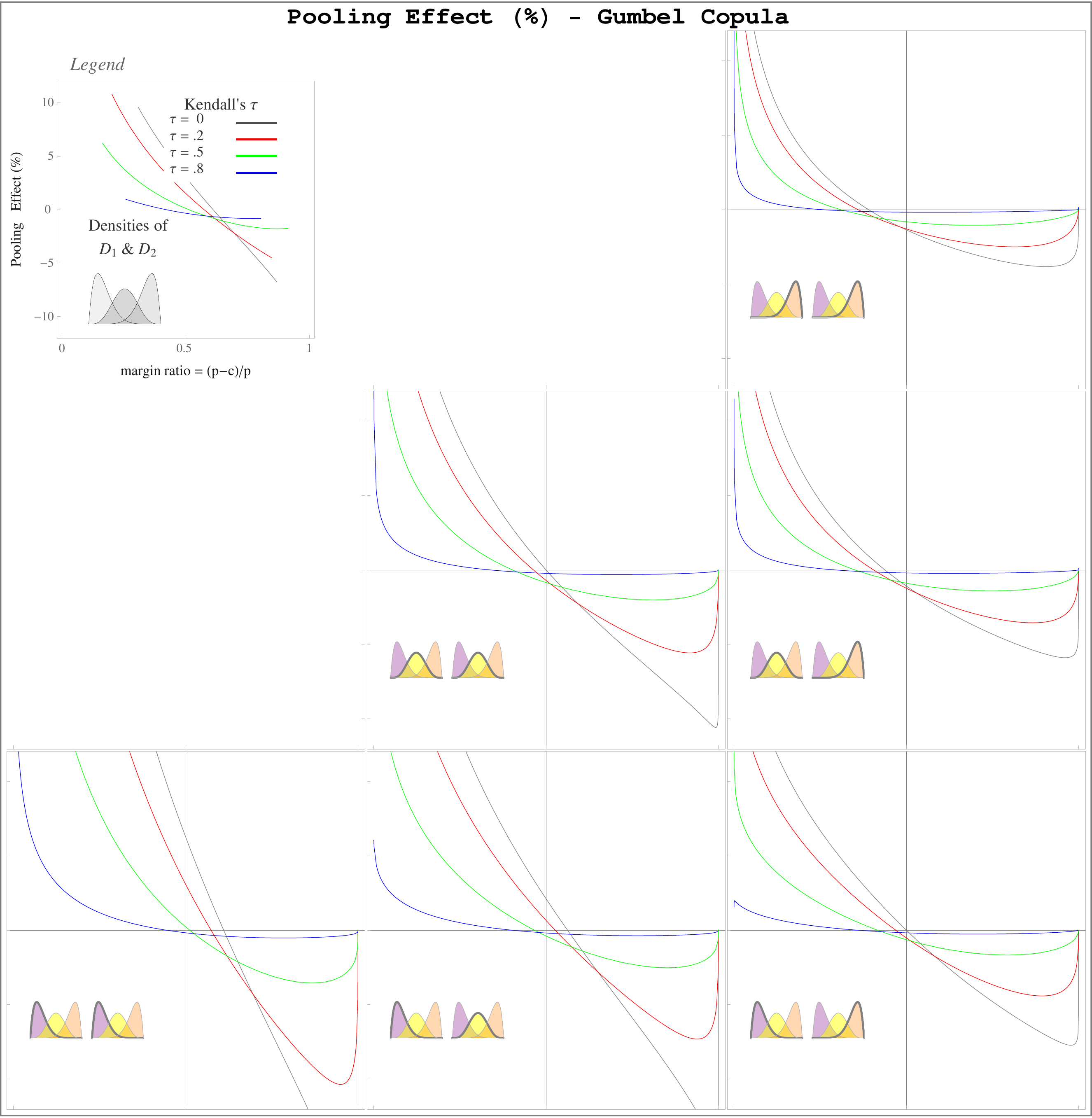}
\end{center}
\caption{The plots of percentage pooling effect on inventory levels under Gumbel copula. The marginal distributions are, from top down, $\beta(2,8)$,$\beta(5,5)$ and $\beta(8,2)$, and from left to right, $\beta(2,8)$,$\beta(5,5)$ and $\beta(8,2)$. These represent left-skewed, symmetrical and right-skewed marginals. The common legend for subplots is given on the top-left corner.}
\end{figure}

\subparagraph{Clayton copula} Compared to Gumbel copula, we find two important differences. First, Clayton implies a higher threshold value $t_0$ compared to Gumbel, given everything the same. Second, this threshold value increases as the dependence increases. This observation is in stark contrast to the one with Gumbel copula. Recall that Clayton copula shifts density towards the left tail of the joint distribution, which is the underlying reason behind these differences.
\begin{figure}
\begin{center}
\includegraphics[scale=0.4]{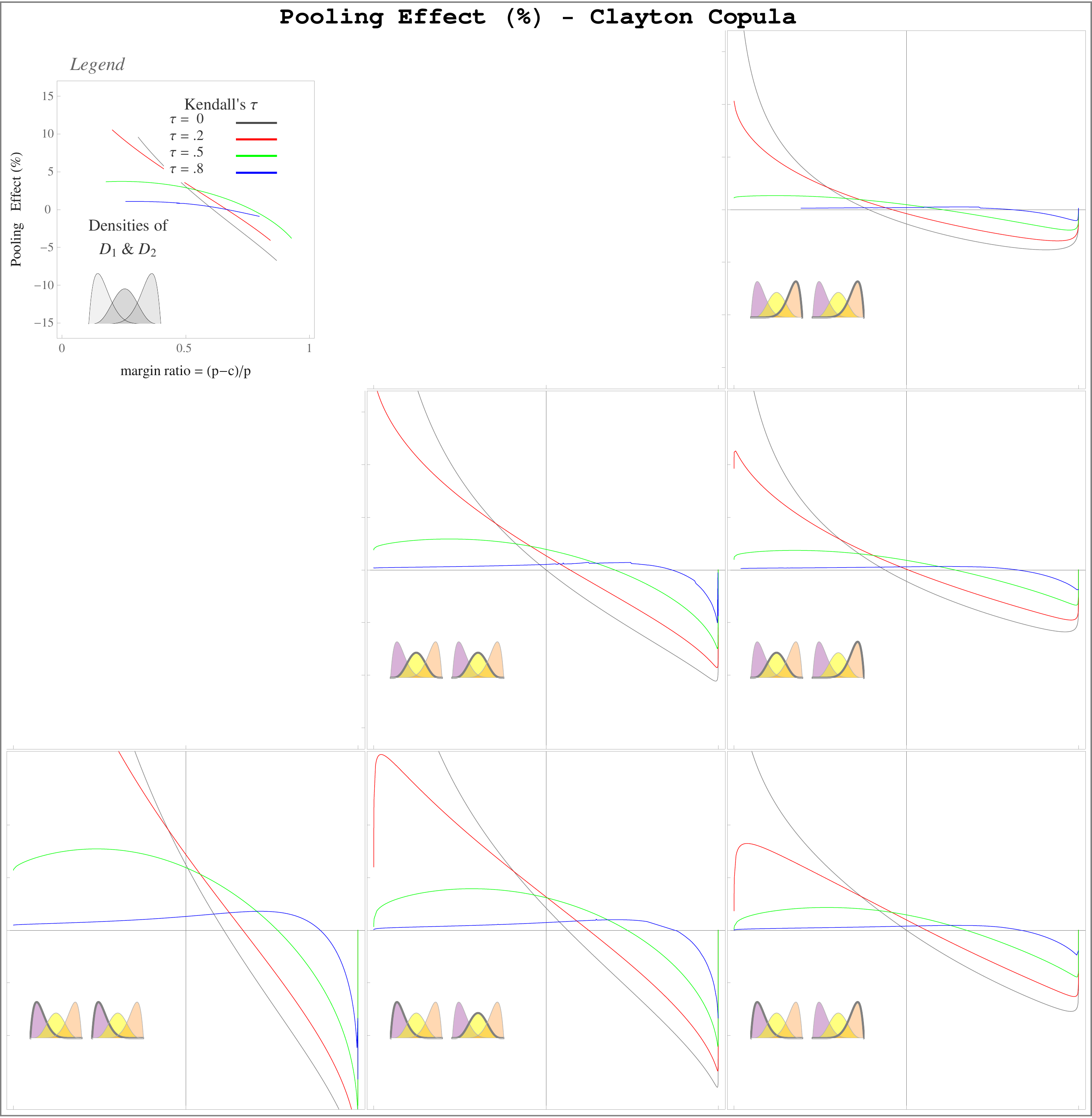}
\end{center}
\caption{The plots of percentage pooling effect on inventory levels under Clayton copula. The marginal distributions are, from top down, $\beta(2,8)$,$\beta(5,5)$ and $\beta(8,2)$, and from left to right, $\beta(2,8)$,$\beta(5,5)$ and $\beta(8,2)$. These represent left-skewed, symmetrical and right-skewed marginals. The common legend for subplots is given on the top-left corner.}
\end{figure}
\subparagraph{Frank copula} This copula can cover both the positive and negative dependence cases. The threshold value $t_0$ is more robust to the skewness of the marginals as compared to the first two copulas. This is due to fact that Frank copula's tail dependencies are symmetric. The effect of this symmetry can also be seen across the pooling effect lines at different $\tau$ levels.
\begin{figure}[ptb]
\begin{center}
\includegraphics[scale=0.4]{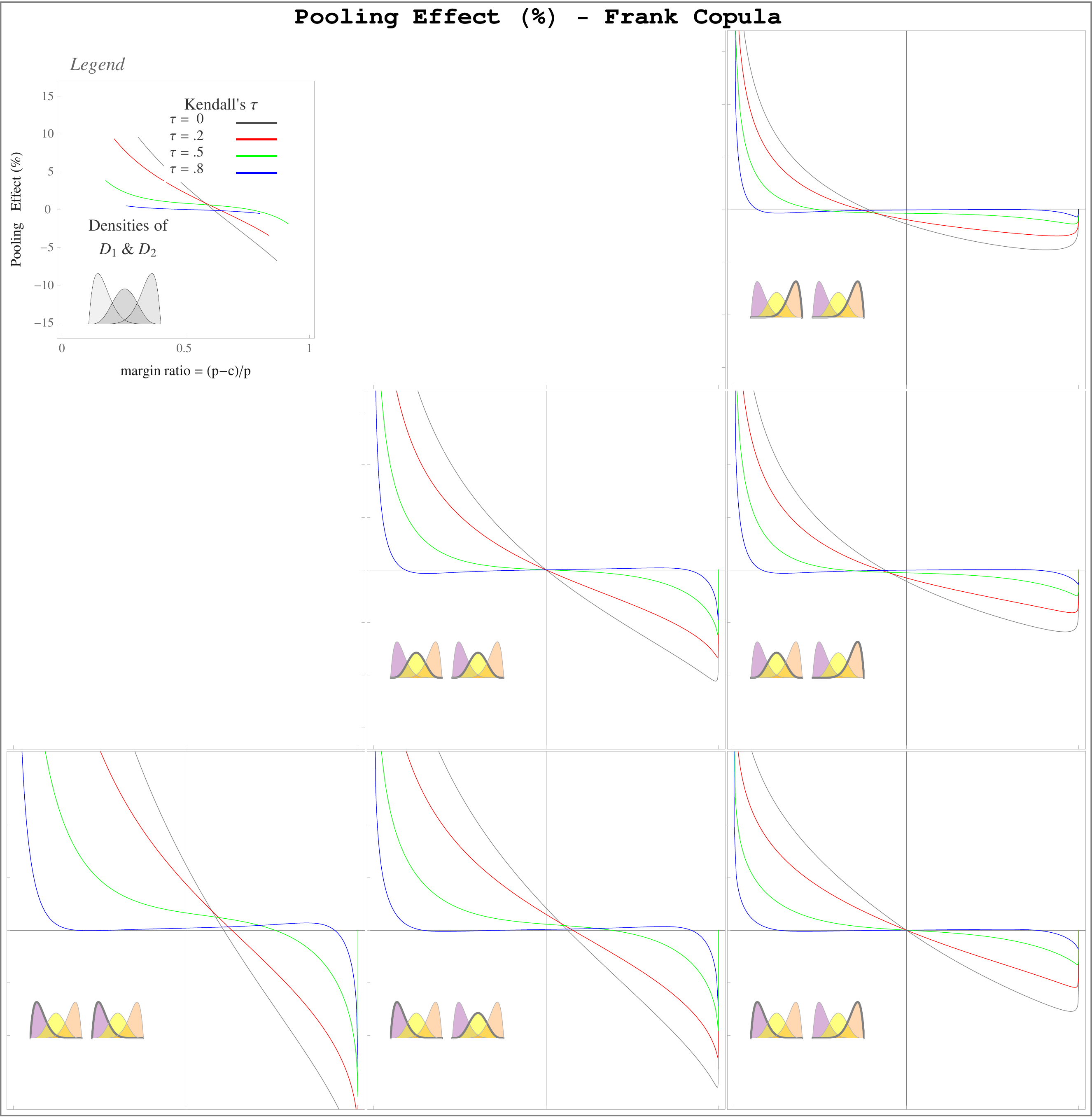}
\end{center}
\caption{The plots of percentage pooling effect on inventory levels under Frank copula with positive Kendall's $\tau$. The marginal distributions are, from top down, $\beta(2,8)$,$\beta(5,5)$ and $\beta(8,2)$, and from left to right, $\beta(2,8)$,$\beta(5,5)$ and $\beta(8,2)$. These represent left-skewed, symmetrical and right-skewed marginals. The common legend for subplots is given on the top-left corner.}
\end{figure}

One case deserves a detailed discussion. When Kendall's $\tau=0.8$, we find that the threshold value is not unique. In Figure 7, we plot the percentile pooling effect for all 6 marginal distribution combinations. On the top of the figure, we plot the regions of the pooling effect where it is positive or negative. For any combination, we find that pooling requires higher inventory levels in two different disjoint regions of the margin ratio. Hence the uniqueness of the threshold value is not valid for this particular copula when there is very high dependence. We should point out that the magnitude of pooling effect is quite small in this case, since co-monotonicity leads to no pooling effect.
\begin{figure}[ptb]
\begin{center}
\includegraphics[scale=0.6]{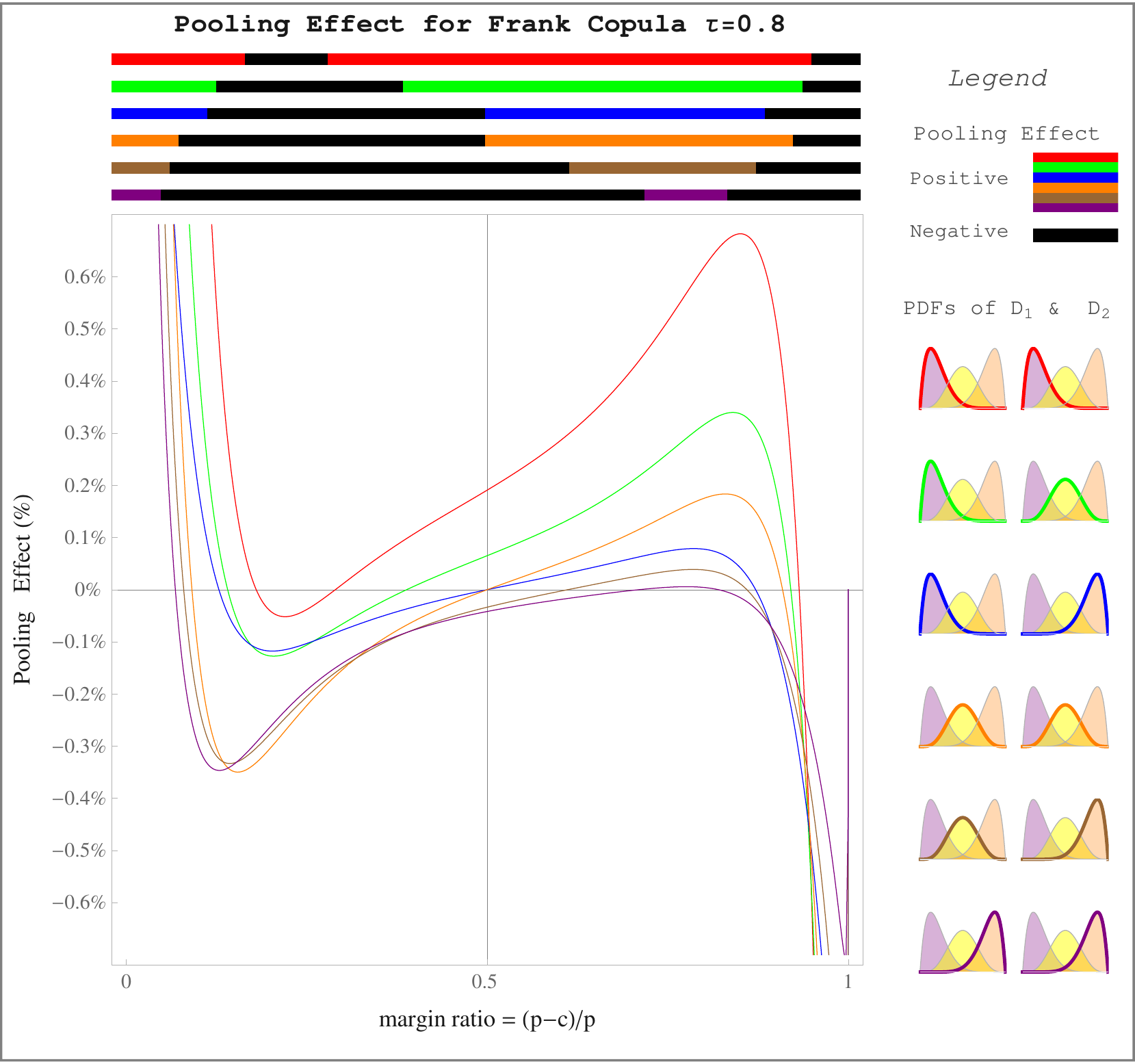}
\end{center}
\caption{The plots of percentage pooling effect on inventory levels under Frank copula with Kendall's $\tau=0.8$. The pooling effect for different marginal distribution combinations, $\beta(2,8)$,$\beta(5,5)$ and $\beta(8,2)$, are depicted in different colors. Top figure shows the regions where pooling effect is negative or positive.}
\end{figure}

With the Frank copula, one can model negative correlation which provides additional insights. As expected, negative correlation can lead to significant pooling effects, especially when the margin ratio is high or low. When demand's are highly negatively correlated, the pooled inventory level (which depends on the sum of these two random variables) is robust against the margin ratio. However, the dedicated inventory levels are small for low margin ratios, and high for high margin ratios. Hence we see that pooling effect is significantly positive for lower ratios and significantly negative for high margin ratios.
\begin{figure}[ptb]
\begin{center}
\includegraphics[scale=0.4]{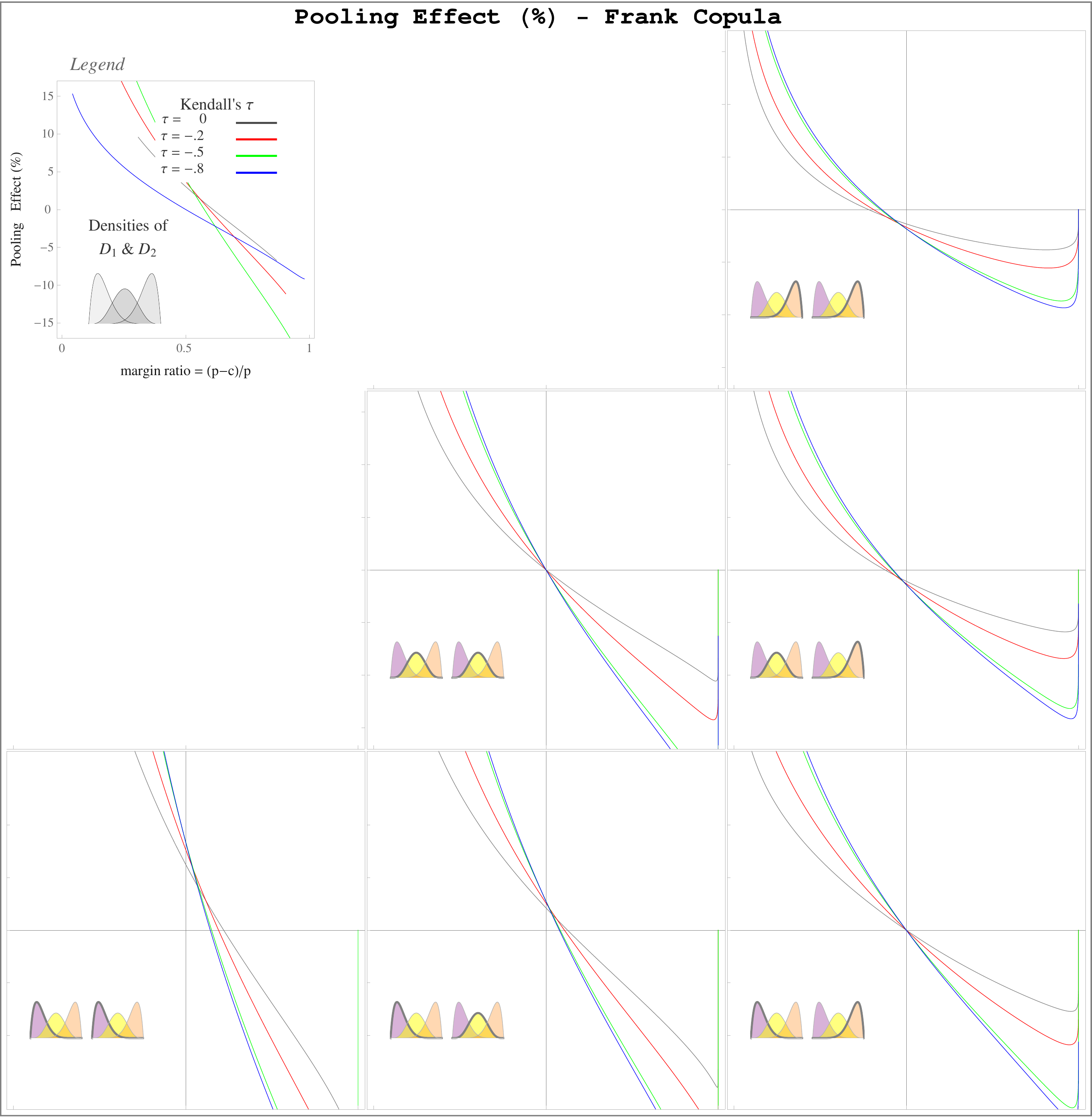}
\end{center}
\caption{The plots of percentage pooling effect on inventory levels under Frank copula with negative Kendall's $\tau$. The marginal distributions are, from top down, $\beta(2,8)$,$\beta(5,5)$ and $\beta(8,2)$, and from left to right, $\beta(2,8)$,$\beta(5,5)$ and $\beta(8,2)$. These represent left-skewed, symmetrical and right-skewed marginals. The common legend for subplots is given on the top-left corner.}
\end{figure}

\subparagraph{Kendall's $\tau$ versus pooling effect} Finally, we address how the magnitude of the pooling effect varies with model parameters. In particular, we investigate whether the optimal pooled inventory level varies monotonically with the dependence.
Under normality, we know that it does. Under general dependence structures, however, monotonicity of pooled inventory level with respect to dependence does not hold. Figure 9 presents the pooled inventory levels under different copulas with identical normal marginals compared to dedicated inventory levels, under three different margin ratios.
\begin{figure}
\begin{center}
\includegraphics[scale=0.5]{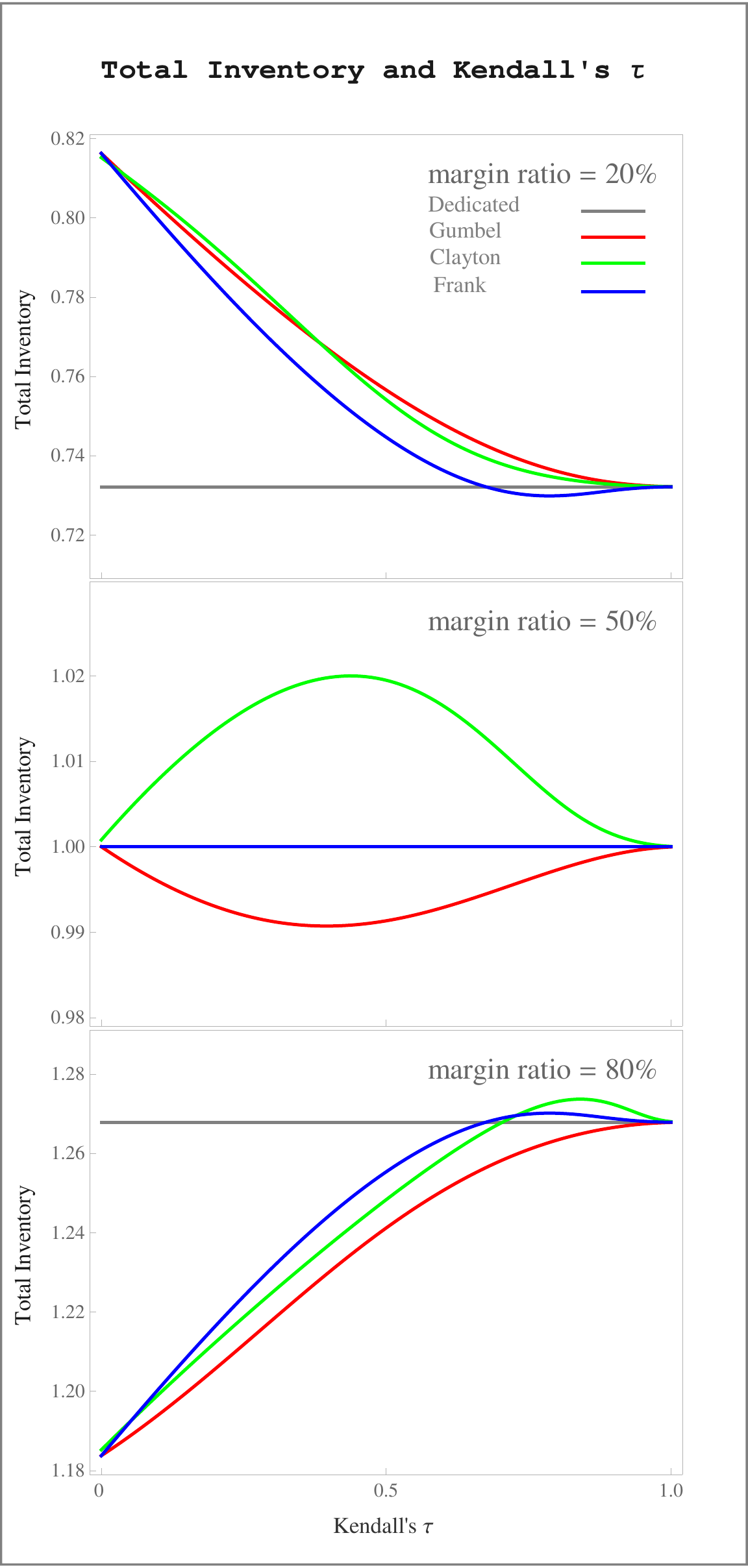}
\caption{The total inventory levels under three different margin ratios, varying with respect to dependence (Kendall's $\tau$). Top panel gives the low margin case with $t=0.2$, middle panel is the medium margin case with $t=0.5$, and the bottom panel is the high margin case with $t=0.8$. The gray line indicates dedicated inventory level, red is the pooled inventory level under Gumbel copula, green is the pooled under Clayton copula and blue line is the pooled inventory under Frank copula. The marginals are $\beta(5,5)$.}
\end{center}
\end{figure}

Figure 9 clearly shows the importance of margin ratio. For a low margin ratio, we see that all $3$ copulas we tried show higher inventory levels compared to dedicated. For high margins, the opposite seems to be the most common trend. For mid-range margin ratios, Clayton copula shows the most variation as the Kendall's $\tau$ changes.

\section{Conclusion}

In this paper, we investigate the optimal inventory levels after pooling to determine whether the manager should increase inventory levels after a switch to pooling. We show that one has to understand the true underlying dependency structure between the individual demand sources, as well as the uncertainty within each demand source, to determine the right level of inventories. In order to understand the effect of each factor, we use copula theory to separate the effect of demand source uncertainty and dependencies between these demand sources, and study the interactions in between, as well as the effects of these interactions on the inventory levels.

We find that the sign of pooling effect depends on the margin ratio: It is positive if and only if margin ratio is higher than a threshold, under certain conditions. This threshold depends on the marginal demand distributions, as well as the copula that joins them. Through numerical studies, we investigate these relationships and conclude that tail dependencies and the strength of dependency are the main factors that affect this threshold. Finally we show that pooled inventory levels are not necessarily monotone with respect to the level of dependence. This is especially true for copulas with asymmetric tail dependencies.

There are open questions that requires further research. In this paper, we focus on the unimodal distributions. If the marginals of the demand distributions are bimodal, then pooling effect is expected to be stronger around these modes, especially if these modes are closer. We also assume that one can estimate the true copula structure underlying the data. When the quality of this estimation is not high or it is not available at all, one will need to consider all possible dependence structures given the limited information, and come up with upper and lower bounds of true optimal inventory levels. Incorrect estimation of the dependence might also cause setting incorrect inventory levels.

In this paper we consider the case with $2$ products. However, the copula framework is able to handle any number of marginals. Therefore is is possible to easily extend the results of this paper to case with arbitrary finite number of products. Our conjecture is that, when the marginals are identical, the results of the multi-item case will be similar to our results when the two marginals are identical. When they are not identical however, the relative shapes of marginals is critical in determining the inventory levels. We leave these questions to future work.

Product characteristics are critical in determining pooled inventory levels. This paper covers the case with perfectly substitutable products. When products are not perfectly substitutable, then pooling effect will be smaller. At the extreme case when pooled products are uniquely different from each other, the effect of pooling should vanish. Hence a deeper understanding of moderate levels of substitutability is required to investigate those situations. Finally, another assumption in this paper is that the products are financially identical: They have the same revenue and cost per unit. When this assumption is not valid, one will need to have a fulfilment policy as to which demand source to satisfy when there is insufficient inventory. This policy structure would determine the pooling inventory levels and in turn the efficiency of pooling.

\begin{acknowledgement}
We would like to thank the two anonymous referees for their helpful comments.
The final version of this manuscript is to appear in Handbook of Newsvendor Problems: Models, Extensions and Applications, Springer, 2012. (www.springerlink.com)
\end{acknowledgement}

\end{document}